\documentclass[a4paper,11pt]{amsart}
\usepackage{fancyhdr}
\usepackage{amsmath}
\usepackage{dsfont}
\usepackage{hyperref}
\usepackage[mathscr]{eucal}
\usepackage[cp1251]{inputenc}
\usepackage[english]{babel}
\usepackage{enumerate,float,indentfirst}
\usepackage{graphicx}
\usepackage{xcolor}
\usepackage{latexsym,a4,mathrsfs,amsthm,amsmath,amssymb,url}
\usepackage{amsfonts}
\usepackage{amssymb}
\usepackage{caption}

\numberwithin{equation}{section}
\newtheorem{formulation}{Formulation}

\begin{document}

\vspace{1in}

\title[Submatrices with the best-bounded inverses]{\bf Submatrices with the best-bounded inverses: revisiting the hypothesis}

\author[Yu. Nesterenko]{Yuri Nesterenko}
\address{ Moscow State University }
\email{y\_nesterenko@mail.ru}

\begin{abstract}

The following hypothesis was put forward by Goreinov, Tyrtyshnikov and Zamarashkin in \cite{GTZ1997}.

For arbitrary real $n \times k$ matrix with orthonormal columns a sufficiently "good" $k \times k$ submatrix exists. "Good" in the sense of having a bounded spectral norm of its inverse. The hypothesis says that for arbitrary $k = 1, \ldots, n-1$ the upper bound can be set at $\sqrt{n}$.

Supported by numerical experiments, the problem remained open for all non-trivial cases ($1 < k < n-1$). In this paper we will give the proof for the simplest of them ($n = 4, \, k = 2$).
\end{abstract}

\maketitle

\thispagestyle{empty}

\section{Introduction}

In this paper we consider the hypothesis proposed by Goreinov, Tyrtyshnikov and Zamarashkin in \cite{GTZ1997}. Related to the theory of pseudoskeleton approximations, initially it was formulated as follows.

\begin{formulation}\label{f1}
For every $n > k > 0$ and an arbitrary real $n \times k$ matrix with orthonormal columns such a $k \times k$ submatrix exists that the spectral norm of its inverse does not exceed $\sqrt{n}$.
\end{formulation}

Let's also give several equivalent formulations\footnote{We leave the proof of equivalence to the reader.}.

\begin{formulation}\label{f2}
For every $n > k > 0$ and an arbitrary real $n \times k$ matrix with orthonormal columns such a $k \times k$ submatrix exists that its $k$-th singular value is not smaller than $1 / \sqrt{n}$.
\end{formulation}

\begin{formulation}\label{f3}
For every $n > k > 0$ and an arbitrary real $k$-dimensional subspace in $\mathds{R}^n$ such a $k$-dimensional coordinate subspace exists that the largest principal angle between them does not exceed $\arccos(1 / \sqrt{n})$.
\end{formulation}

\begin{formulation}\label{f4}
For every $n > k > 0$ and an arbitrary $k$-circle centered at the origin of $\mathds{R}^n$, there exists an orthogonal projection onto a coordinate $k$-subspace lying entirely outside the origin centered open ball of radius $1 / \sqrt{n}$.
\end{formulation}

Our numerical experiments support this hypothesis and, moreover, indicate that the provided estimate is sharp. In the next section we will prove this fact for $n = 4, \, k = 2$.

In what follows, we will use Formulation \ref{f2}. 

\section{Proof for $n = 4, \, k = 2$}

Let's consider matrix $A \in \mathds{R}^{4 \times 2}$ such that
\begin{equation}\label{eq:0}
A^T A = I.
\end{equation}

Its Pl\"ucker coordinates $p_{ij}$ computed as the minors in the corresponding rows satisfy the Pl\"ucker relation\footnote{In our case it can be simply derived from the equation $\det(A|A) = 0$.}
\begin{equation}\label{eq:033}
p_{12} p_{34} - p_{13} p_{24} + p_{14} p_{23} = 0,
\end{equation}
as well as the normalization condition
\begin{equation}\label{eq:066}
p_{12}^2 + p_{13}^2 + p_{14}^2 + p_{23}^2 + p_{24}^2 + p_{34}^2 = 1,
\end{equation}
following from the equation (\ref{eq:0}) and the Binet-Cauchy identity.

Let's find out what values the Pl\"ucker coordinates can take if the minimum singular values of all $2 \times 2$ submatrices of $A$ are bounded above by the value $1/2$.

To do this, we use the thin version of $CS$-decomposition (see \cite{Golub2013})
\begin{equation*}
A = \left( \begin{matrix}
Q_1 & \quad \\
\quad & Q_2
\end{matrix} \right)
\left( \begin{matrix}
c_1 & 0 \\
0 & c_2 \\
s_1 & 0 \\
0 & s_2 \\
\end{matrix} \right) Q_3.
\end{equation*}
Here, $c_1, c_2$ and $s_1, s_2$ are the corresponding cosines and sines of some arguments $0 \leq \alpha \leq \beta \leq \pi/2$, and $Q_1, Q_2, Q_3$ are orthogonal $2 \times 2$ matrices.

From the considered upper bound for the minimum singular value, we derive that $|p_{12}|$ and $|p_{34}|$ can take the values $\cos\alpha \cos\beta$ and $\sin\alpha \sin\beta$, where $0 \leq \alpha \leq \pi/6$ and $\pi/3 \leq \beta \leq \pi/2$. Elimination of parameters $\alpha$ and $\beta$ yields the system of inequalities
\begin{equation*}
\left\{ \begin{split}
&4 p_{12}^2 + \frac{4}{3} p_{34}^2 \leq 1, \\
&\frac{4}{3} p_{12}^2 + 4 p_{34}^2 \leq 1.
\end{split}\right.
\end{equation*}
Similar inequalities are also valid for the remaining pairs $p_{13}, p_{24}$ and $p_{14}, p_{23}$.

Let's rewrite the resulting system of six inequalities together with the equations (\ref{eq:033}) and (\ref{eq:066}) in terms of the variables
\begin{equation*}
\begin{split}
&x_1 = p_{12} + p_{34}, \quad x_2 = p_{12} - p_{34}, \\
&y_1 = p_{13} - p_{24}, \quad y_2 = p_{13} + p_{24}, \\
&z_1 = p_{14} + p_{23}, \quad z_2 = p_{14} - p_{23},
\end{split}
\end{equation*}
(the so-called Klein coordinates multiplied by 2). After a simple transformation it yields
\begin{equation}\label{eq:1}
\left\{ \begin{split}
&x_1^2 + y_1^2 + z_1^2 = 1, \\
&x_2^2 + y_2^2 + z_2^2 = 1, \\
&x_1^2 \pm x_1 x_2 + x_2^2 \leq 3/2, \\
&y_1^2 \pm y_1 y_2 + y_2^2 \leq 3/2, \\
&z_1^2 \pm z_1 z_2 + z_2^2 \leq 3/2.
\end{split}\right.
\end{equation}

Since the system (\ref{eq:1}) is invariant under change of sign of its unknowns, it is sufficient to examine its consistency in the case of non-negative solutions.

In this case, we can apply the following change of variables
\begin{equation*}
\left\{ \begin{split}
&x_1 = X \sin(x + \pi/3), \quad x_2 = X \sin(x - \pi/3), \\
&y_1 = Y \sin(y + \pi/3), \quad y_2 = Y \sin(y - \pi/3), \\
&z_1 = Z \sin(z + \pi/3), \quad z_2 = Z \sin(z - \pi/3),
\end{split}\right.
\end{equation*}
parametrizing the elliptic sectors given by (\ref{eq:1}). As a result, the corresponding inequalities become much simpler
\begin{equation*}
\left\{ \begin{split}
&0 \leq X, Y, Z \leq 1, \\
&\frac{\pi}{3} \leq x, y, z \leq \frac{2\pi}{3},
\end{split}\right.
\end{equation*}
while the equations take the form
\begin{equation}\label{eq:2}
\left\{ \begin{split}
&X^2 \sin^2(x + \pi/3) + Y^2 \sin^2(y + \pi/3) + Z^2 \sin^2(z + \pi/3) = 1, \\
&X^2 \sin^2(x - \pi/3) + Y^2 \sin^2(y - \pi/3) + Z^2 \sin^2(z - \pi/3) = 1.
\end{split}\right.
\end{equation}

Let us show that the given system of equations and inequalities is consistent if and only if $X = Y = Z = 1$. This will prove the estimate we are interested in this paper, as well as its sharpness.

From the equations \eqref{eq:2} and the inequalities $X, Y, Z \leq 1$ we have
\begin{equation}\label{eq:3}
\left\{ \begin{split}
&\sin^2(x + \pi/3) + \sin^2(y + \pi/3) + \sin^2(z + \pi/3) \geq 1, \\
&\sin^2(x - \pi/3) + \sin^2(y - \pi/3) + \sin^2(z - \pi/3) \geq 1.
\end{split}\right.
\end{equation}

Under the existing constraints $\frac{\pi}{3} \leq x, y, z \leq \frac{2\pi}{3}$, the first of the conditions \eqref{eq:3} leads to the inequality $x + y + z \leq 3\pi/2$,
and the second one --- to the opposite inequality $x + y + z \geq 3\pi/2$. Below we prove the second\footnote{The first implication can be proved in a similar way.} of these implications by assuming the contrary.

From the given assumption and from the monotonicity of sine function on the segment $[0, \pi/3]$, there exist values $x', y', z' \in [0, \pi/3]$ satisfying both the equality $x' + y' + z' = \pi/2$ and the inequality
\begin{equation*}
\sin^2x' + \sin^2y' + \sin^2z' > 1.
\end{equation*}

However, analysis of the function $\sin^2x' + \sin^2y' + \sin^2z'$ inside the triangle $x' + y' + z' = \pi/2, \, 0 \leq x', y', z' \leq \pi/2$ shows that it attains its maximum on the boundary, where it is identically equal to $1$. Therefore, we have got a contradiction.

Now, the assumption that at least one of the parameters $X, Y, Z$ is strictly less than $1$ makes at least one of the inequalities \eqref{eq:3} strict,
which makes this system inconsistent. The whole proof is complete.

\begin{figure}[H]
\centering
\includegraphics[width=12cm]{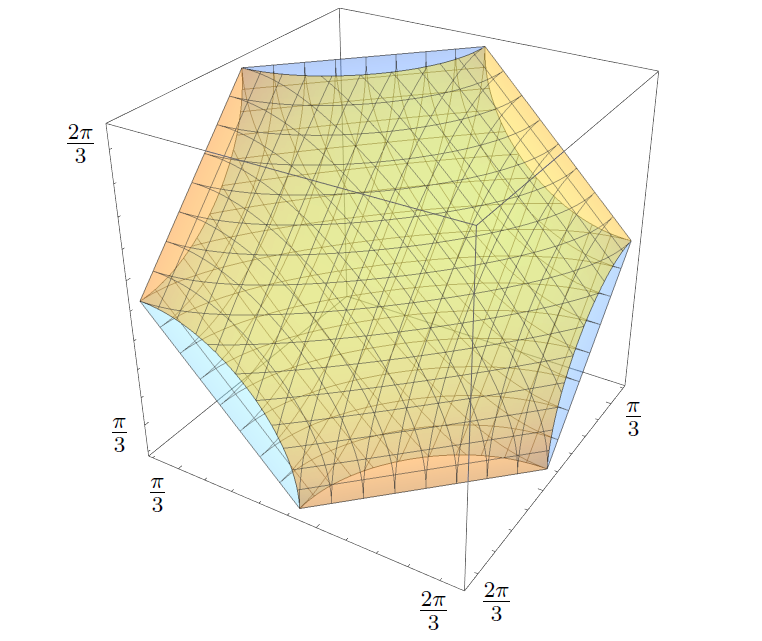}
\captionsetup{justification=centering}
\caption{ The boundaries of the regions \eqref{eq:3}. The contact points correspond to the subspaces the most deviating from the coordinate ones. }
\label{fig:plot}
\end{figure}
\vspace*{-0.25cm}
Solving the system \eqref{eq:1} for $X = Y = Z = 1$, it is possible to construct matrices for which the estimate considered in the hypothesis turns into equality. All the~corresponding $96$ subspaces (in terms of Formulation \ref{f3}) may be obtained from the one given by the matrix
\begin{equation*}
\left( \begin{matrix}
\sqrt{\frac{1}{2}} &\sqrt{\frac{1}{8}} \\
-\sqrt{\frac{1}{2}} &\sqrt{\frac{1}{8}} \\
0 &\sqrt{\frac{3}{8}} \\
0 &\sqrt{\frac{3}{8}}
\end{matrix} \right)
\end{equation*}
using symmetry transformations.

\section{Acknowledgements}

I would like to thank the authors of the original paper \cite{GTZ1997} for a lot of fruitful discussions as well as Igor Makhlin who proposed Formulation~\ref{f4}.

\bibliographystyle{plain}
\bibliography{lit}

\end{document}